\documentclass[11pt]{article}

\usepackage{amsmath}
\usepackage{amssymb}
\usepackage{indentfirst}
\usepackage{graphics} 
\usepackage{color}

\makeatletter

\@addtoreset{equation}{section}
\makeatother
\def\<{\langle}
\def\>{\rangle}

\newtheorem{lem}{Lemma}[section]
\newtheorem{theo}{Theorem}[section]
\newtheorem{rem}{Remark}[section]
\newtheorem{pro}{Proposition}[section]

\makeatletter
   
   \@addtoreset{equation}{section}
\makeatother

\begin{document}
\title{\bf Asymptotic profile and optimal decay of solutions \\ of some wave equations with logarithmic damping}
\author{Ruy Coimbra Char\~ao\thanks{ruy.charao@ufsc.br}  \\{\small Department of Mathematics} \\{\small Federal University of Santa Catarina} \\ {\small 88040-270, Florianopolis, Brazil} \\and\\Ryo Ikehata\thanks{Corresponding author: ikehatar@hiroshima-u.ac.jp} \\ {\small Department of Mathematics}\\ {\small Graduate School of Education} \\ {\small Hiroshima University} \\ {\small Higashi-Hiroshima 739-8524, Japan}}
\date{}
\maketitle
\begin{abstract}
We introduce a new model of the nonlocal wave equation with a logarithmic damping mechanism, which is rather weak as compared with frequently studied fractional damping cases.  We consider the Cauchy problem for the new model in ${\bf R}^{n}$. We study the asymptotic profile and optimal decay rates of solutions as $t \to \infty$ in $L^{2}$-sense. The damping terms considered in this paper is not studied so far, and in the low frequency parameters the damping is rather weakly effective than that of well-studied power type one such as $(-\Delta)^{\theta}u_{t}$ with $\theta \in (0,1)$. In order to get the optimal rate of decay we meet the so-called hypergeometric functions, so the analysis seems to be more difficult and attractive.  
\end{abstract}
\section{Introduction}
\footnote[0]{Keywords and Phrases: Wave equation; Logarithmic damping; $L^{2}$-decay; asymptotic profile, optimal decay.}
\footnote[0]{2010 Mathematics Subject Classification. Primary 35L05; Secondary 35B40, 35C20, 35S05.}
We present and consider a new type of wave equation with a logarithmic damping term:
\begin{equation}
u_{tt} + Au + Lu_{t} = 0,\ \ \ (t,x)\in (0,\infty)\times {\bf R}^{n},\label{eqn}
\end{equation}
\begin{equation}
u(0,x)= u_{0}(x),\ \ u_{t}(0,x)= u_{1}(x),\ \ \ x\in{\bf R}^{n} ,\label{initial}
\end{equation}
where $(u_{0},u_{1})$ are initial data chosen as
\[u_{0} \in H^{1}({\bf R}^{n}),\quad u_{1} \in L^{2}({\bf R}^{n}),\]
and the operator $Au := -\Delta u$ for $u \in H^{2}({\bf R}^{n})$, and a new operator 
\[L: D(L) \subset L^{2}({\bf R}^{n}) \to L^{2}({\bf R}^{n})\]
is defined as follows: 
\[D(L) := \left\{f \in L^{2}({\bf R}^{n}) \,\bigm|\,\int_{{\bf R}^{n}}(\log(1+\vert\xi\vert^{2}))^{2}\vert\hat{f}(\xi)\vert^{2}d\xi < +\infty\right\},\]
for $f \in D(L)$,  
\[(Lf) (x) := {\cal F}_{\xi\to x}^{-1}\left(\log (1+\vert\xi\vert^{2})\hat{f}(\xi)\right)(x),\]
and symbolically writing, one can see
\[L = \log(I+A).\] 

Here, we denote the Fourier transform ${\cal F}_{x\to\xi}(f)(\xi)$ of $f(x)$ by 
\[{\cal F}_{x\to\xi}(f)(\xi) = \hat{f}(\xi) := \displaystyle{\int_{{\bf R}^{n}}}e^{-ix\cdot\xi}f(x)dx\]
as usual with $i := \sqrt{-1}$, and ${\cal F}_{\xi\to x}^{-1}$ expresses its inverse Fourier transform. Since the new operator $L$ is constructed by a nonnegative-valued multiplication one, it is nonnegative and self-adjoint in $L^{2}({\bf R}^{n})$. 

Then, by a similar argument to \cite[Proposition 2.1]{ITY} based on Lumer-Phillips Theorem one can find that the problem (1.1)-(1.2) has a unique mild solution
\[u \in C([0,\infty);H^{1}({\bf R}^{n})) \cap C^{1}([0,\infty);L^{2}({\bf R}^{n}))\]
satisfying the energy inequality
\begin{equation}\label{energy}
E_{u}(t) \leq E_{u}(0),
\end{equation} 
where
\[
E_{u}(t) := \frac{1}{2}\left(\Vert u_{t}(t,\cdot)\Vert_{L^{2}}^{2} + \Vert\nabla u(t,\cdot)\Vert_{L^{2}}^{2}\right).
\]
\eqref{energy} implies the decreasing property of the total energy because of the existence of some kind of dissipative term $Lu_{t}$. For details, see Appendix in this paper.

A main topic of this paper is to find an asymptotic profile of solutions in the $L^{2}$ topology as $t \to \infty$ to problem (1.1)-(1.2), and to apply it to get the optimal rate of decay of solutions in terms of the $L^{2}$-norm.\\  

Now let us recall several previous works related to linear damped wave equations with constant coefficients. We mention them from the viewpoint of the Fourier transformed equations.\\

(1)\,In the case of weak damping for the equation (1.1) with $L = I$:
\begin{equation}\label{weak}
\hat{u}_{tt} + \vert\xi\vert^{2}\hat{u} + \hat{u}_{t} = 0,
\end{equation}  
as is usually observed, for small $\vert\xi\vert$, the solution to the equation \eqref{weak} behaves like a constant multiple of the Gauss kernel, while for large $\vert\xi\vert$ the solution includes an oscillation property, which vanishes very fast. This property can be pointed out in the celebrated paper due to Matsumura \cite{M} from the precise decay estimates point of view, and from the viewpoint of asymptotic expansions one should mention the so-called Nishihara decomposition \cite{Ni}. As a work in a similar philosophy one can cite the paper due to Narazaki \cite{Na}, which deals with the higher dimensional case. Higher order asymptotic expansions in $t$ of the $L^{2}$-norm of solutions to \eqref{weak} has been studied in Volkmer \cite{V}, and Said-Houari \cite{SH} restudies the diffusion phenomena from the viewpoint of the weighted $L^{1}$-initial data. As an abstract theory in the case when the operator $A$ is nonnegative and self-adjoint in Hilbert spaces one can cite several papers due to Chill-Haraux \cite{CH}, Ikehata-Nishihara \cite{IN}, Radu-Todorova-Yordanov \cite{RTY}, and Sobajima \cite{SJ}, where they also investigate the diffusion phenomenon of solutions (as $t \to \infty$) more precisely.\\

(2)\,After the equation \eqref{weak} one should give some comments to the following Fourier transformed equation of (1.1) with $L = A$:          
\begin{equation}\label{strong}
\hat{u}_{tt} + \vert\xi\vert^{2}\hat{u} + \vert\xi\vert^{2}\hat{u}_{t} = 0.
\end{equation} 
This corresponds to the so-called strongly damped wave equation case, which was set up by Ponce \cite{Po} and Shibata \cite{S} at the first stage of the research. On the contrary to \eqref{weak} the solution to the equation \eqref{strong} behaves like a diffusion wave with complex-valued characteristic roots for small $\vert \xi\vert$, while for large $\vert\xi\vert$ the solution does not include any oscillation property because the characteristic roots are real-valued. These observations are recently pointed out in the papers Ikehata \cite{I-14}, Ikehata-Onodera \cite{IO-17} and Ikehata-Todorova-Yordanov \cite{ITY} by capturing the leading term as $t \to \infty$ of the solution, and they have derived optimal estimates of solutions in terms of $L^{2}$-norm. In particular, it should be noted that optimal estimates in \cite{I-14, IO-17} do not necessarily imply the decay estimates, and in fact, in space dimension $1$ and $2$ an infinite time blowup property occurs. Furthermore, it should be emphasized that higher order asymptotic expansions of the squared $L^{2}$-norm of solutions as $t \to \infty$ have been very precisely studied recently by Barrera \cite{B} and Barrera-Volkmer \cite{BV-1, BV-2}, and in \cite{Mi} Michihisa studies the higher order asymptotic expansion of the solution itself, and applied it to investigate the lower bound of decay rate of the difference between the leading term and the solution itself in terms of $L^{2}$ norm. \\

(3)\,As for the generalization of the study (1) and (2) one can cite several papers due to D'Abbicco-Ebert \cite{DE}, D'Abbicco-Ebert-Picon \cite{DEP}, D'Abbicco-Reissig \cite{DR}, Char\~ao-da Luz-Ikehata \cite{CLI}, Ikehata-Takeda \cite{IT}, Karch \cite{K} and Narazaki-Reissig \cite{NR}, which deal with the equation (1.1) with the so-called structural damping $L = A^{\theta}$ ($0 < \theta < 1$):    
\begin{equation}\label{structural}
\hat{u}_{tt} + \vert\xi\vert^{2}\hat{u} + \vert\xi\vert^{2\theta}\hat{u}_{t} = 0.
\end{equation} 
By those works one has already known that $\theta = 1/2$ is critical in the sense that for $\theta \in (0,1/2)$ the solution to the equation \eqref{structural} is parabolic like, and in the case of $\theta \in [1/2, 1)$ the solution to the equation \eqref{structural} behaves like a diffusion wave as $t \to \infty$, and in particular, $\theta = 1/2$ corresponds to the scale invariant case. \\

On reconsidering our problem in the Fourier space our equation (1.1) becomes  
\begin{equation}\label{logdamping}
\hat{u}_{tt} + \vert\xi\vert^{2}\hat{u} + \log(1+\vert\xi\vert^{2})\hat{u}_{t} = 0.
\end{equation} 
An influence of the damping coefficient $\log(1+\vert\xi\vert^{2})$ on the dissipative nature of the solution seems to be rather weak in both regions for large $\vert\xi\vert$ and small $\vert\xi\vert$ when we compare it with the fractional damping $\vert\xi\vert^{2\theta}$ with $\theta \in (1/2,1)$. Furthermore, as is easily seen that the characteristic roots $\lambda_{\pm}$ for the characteristic polynomial of \eqref{logdamping} such that 
\[\lambda^{2} + \log(1+\vert\xi\vert^{2})\lambda + \vert\xi\vert^{2} = 0\]
are all complex-valued for all $\xi \in {\bf R}^{n}$, which implies that the corresponding solution has an oscillating property for all frequency parameter $\xi \in {\bf R}^{n}$. This property produces a big difference as compared with previously considered cases in the references. So, the problem (1.1) may include several difficulties which has never experienced so far, and in fact, we have to analyze the so-called hypergeometric functions naturally (see section 2) when one captures the leading term and obtains the optimal rate of decay of the solution to problem (1.1)-(1.2). In this sense, we do present much new type of problems in the wave equation field through the analysis for the wave equation with $\log$-damping. \\ 

Our two main results read as follows.

\begin{theo}\label{main-theo}
\, Let $n \geq 1$, and let $[u_{0},u_{1}] \in \left(H^{1}({\bf R}^{n})\cap L^{1}({\bf R}^{n})\right) \times \left(L^{2}({\bf R}^{n})\cap L^{1,1}({\bf R}^{n})\right)$. Then, the unique solution $u(t,x)$ to problem {\rm (1.1)}-{\rm (1.2)} satisfies
\[\left\Vert u(t,\cdot) - \left(\int_{{\bf R}^{n}}u_{1}(x)dx\right){\cal F}_{\xi\to x}^{-1}\left((1+\vert\xi\vert^{2})^{-\frac{t}{2}}\frac{\sin(\vert\xi\vert t)}{\vert\xi\vert}\right)\right\Vert_{L^{2}} \leq I_{0}t^{-\frac{n}{4}}, \quad (t \gg 1),\]
where
\[I_{0} := \Vert u_{0}\Vert_{L^{2}} + \Vert u_{1}\Vert_{L^{2}} + \Vert u_{0}\Vert_{L^{1}} + \Vert (1+\vert x\vert)u_{1}\Vert_{L^{1}}.\]
\end{theo}
\begin{rem}\,{\rm $I_{0}$ in Theorem \ref{main-theo} does not depend on any norms $\Vert u_{0}\Vert_{H^{1}}$. This is one of differences as compared with the case for fractional damping $L := (-\Delta)^{\theta}$ (see \cite{I-14}). }
\end{rem}

As a consequence of Theorem \ref{main-theo} one can get the optimal estimates in $t$ of solutions in terms of $L^{2}$-norm. We set
\[P_{1} := \int_{{\bf R}^{n}}u_{1}(x)dx.\]
\begin{theo}\label{main-theo2}
\, Let $n \geq 1$, and let $[u_{0},u_{1}] \in \left(H^{1}({\bf R}^{n})\cap L^{1}({\bf R}^{n})\right) \times \left(L^{2}({\bf R}^{n})\cap L^{1,1}({\bf R}^{n})\right)$. Then, the unique solution $u(t,x)$ to problem {\rm (1.1)}-{\rm (1.2)} satisfies\\
\vspace{0.1cm}
{\rm (i)}\,\,\,$n \geq 3$ $\Rightarrow$ $C_{n}\vert P_{1}\vert t^{-\frac{n-2}{4}} \leq \Vert u(t,\cdot)\Vert_{L^{2}} \leq C_{n}^{-1}I_{0}t^{-\frac{n-2}{4}}$ {\rm (}$t \gg 1${\rm )},\\
\vspace{0.1cm}
{\rm (ii)}\,\,$n = 2$ $\Rightarrow$ $C_{2}\vert P_{1}\vert \sqrt{\log t} \leq \Vert u(t,\cdot)\Vert_{L^{2}} \leq C_{2}^{-1}I_{0}\sqrt{\log t}$ {\rm (}$t \gg 1${\rm )},\\
\vspace{0.1cm}
{\rm (iii)}\,$n = 1$ $\Rightarrow$ $C_{1}\vert P_{1}\vert \sqrt{t} \leq \Vert u(t,\cdot)\Vert_{L^{2}} \leq C_{1}^{-1}I_{0}\sqrt{t}$ {\rm (}$t \gg 1${\rm )},\\
where $I_{0}$ is a constant defined in Theorem {\rm \ref{main-theo}}, and $C_{n}$ {\rm (}$n \in {\bf N}${\rm )} are constants independent from any $t$ and initial data.  
\end{theo}
\begin{rem}\,{\rm In the case of $\vert\xi\vert \approx 0$ one has
\[\log(1+\vert\xi\vert^{2}) \leq \vert\xi\vert^{2\theta}\]
for all $\theta \in (0,1)$. This implies, in the case of $\theta \in (0,1]$, the effect of damping $\log(1+\vert\xi\vert^{2})\hat{u}_{t}$ is much weaker than all types of fractional damping $\vert\xi\vert^{2\theta}\hat{u}_{t}$. While, if we compare two types of fractional damping $(-\Delta)^{\theta}w_{t}$ with $\theta \in (0,1/2)$ (parabolic like) and $(-\Delta)^{\theta}w_{t}$ with $\theta \in (1/2,1]$ (diffusion wave like), in the case of small $\vert\xi\vert$ (this is essential part of both solutions), the effect of $(-\Delta)^{\theta}w_{t}$ with $\theta \in (1/2,1)$ is weaker than $(-\Delta)^{\theta}w_{t}$ with $\theta \in (0,1/2)$. In some sense, since the effect of damping $\log(1+\vert\xi\vert^{2})\hat{u}_{t}$ is weaker than $\vert\xi\vert^{2\theta}\hat{u}_{t}$ with $\theta \in (1/2,1]$ (wave like case), it seems to be natural that the results obtained in Theorem \ref{main-theo2} coincide with the case for strong damping $(-\Delta)u_{t}$ (see \cite{I-14, IO-17}), however, its analysis is much more difficult, in particular, when we deal with the upper bound for several quantities. This is because the ingredient $(1+\vert\xi\vert^{2})^{-\frac{t}{2}}$ of the leading term obtained in Theorem \ref{main-theo} behaves slower than usual diffusion wave case such that $e^{-t\vert\xi\vert^{2}}\frac{\sin(\vert\xi\vert t)}{\vert\xi\vert}$ in the Fourier space ${\bf R}_{\xi}^{n}$ for each $t > 0$. As will be observed in section 2, the function $(1+\vert\xi\vert^{2})^{-\frac{t}{2}} = e^{-\frac{t}{2}\log(1+\vert\xi\vert^{2})}$ has a close relation to the so-called hypergeometric function with special parameters. In this connection, it would be interesting to study a kind of diffusion equation such that
\[v_{t} + \frac{1}{2}\log(I-\Delta)v = 0,\]
in order to know more about deeper properties of the asymptotic profile obtained in Theorem \ref{main-theo}.}
\end{rem}
\par
\vspace{0.1cm}

This paper is organized as follows. In section 2 we prepare several important propositions and lemmas, which will be used later, and in particular, in subsection 2.1 we shall mention the so-called hypergeometric functions. Theorem \ref{main-theo} is proved in section 3. In section 4, we shall study the optimality of the $L^{2}$-norm of solutions to problem (1.1)-(1.2) in the case of space dimension $1$ and $2$, and Theorem \ref{main-theo2} will be proved at a stroke. Appendix is prepared to check the unique existence of the weak solution to problem (1.1)-(1.2).\\

{\bf Notation.} {\small Throughout this paper, $\| \cdot\|_q$ stands for the usual $L^q({\bf R}^{n})$-norm. For simplicity of notation, in particular, we use $\| \cdot\|$ instead of $\| \cdot\|_2$. Furthermore, we denote $\Vert\cdot\Vert_{H^{l}}$ as the usual $H^{l}$-norm. Furthermore, we define a relation $f(t) \sim g(t)$ as $t \to \infty$ by: there exist constant $C_{j} > 0$ ($j = 1,2$) such that
\[C_{1}g(t) \leq f(t) \leq C_{2}g(t)\quad (t \gg 1).\] 

We also introduce the following weighted functional spaces.
\[L^{1,\gamma}({\bf R}^{n}) := \left\{f \in L^{1}({\bf R}^{n}) \; \bigm| \; \Vert f\Vert_{1,\gamma} := \int_{{\bf R}^{n}}(1+\vert x\vert^{\gamma})\vert f(x)\vert dx < +\infty\right\}.\]
Finally, we denote the surface area of the $n$-dimensional unit ball by $\omega_{n} := \displaystyle{\int_{\vert\omega\vert = 1}}d\omega$. 

}


\section{Hypergeometric functions}

Our interest refers to the historically most important hypergeometric function ${}_2F_1(a,b;c;z)$ called Gauss's hypergeometric function which may be defined by
\begin{equation}\label{hypergeomet}
{ }_{2}F_1(a,b;c;z)=\sum_{n=0}^{\infty} \frac{(a)_n(b)_n}{(c)_n}\frac{z^n}{n!},
\end{equation}
where $(a)_n=a(a+1) \cdots (a+n-1)$ is the Pochhammer symbol (upward factorial). The series \eqref{hypergeomet} converges absolutely in $\vert z\vert < 1$ ($z \in {\bf C}$) for parameters $a,b,c \in {\bf C}$ with $c\neq 0, -1, -2,-3, \cdots$

The generalized hypergeometric functions ${}_pF_q(a_1, \cdots,a_p, b_1, \cdots, b_p; c_1,\cdots,c_q;z)$ are also defined  similarly by hypergeometric power series, that include many other special functions as, for example,  Beta function.

These functions appear in many problems in statistics, probability, quantum mechanics among other areas.
For a list of some of the many thousands of published identities, symmetries, limits, involving the hypergeometric functions we can refer to the works by Erd\'elyi et al. (1955),  Gasper- Rahman (2004), Miller-Paris (2011). There is no any known systems for organizing all of the identities. In fact, there is no known algorithm that can generate all identities. Moreover, there are known a number of different algorithms that generate different series of identities. The theory about the  algorithmic remains an active research topic. 

Hypergeometric function is also given as a solution of the special Euler second-order linear ordinary differential equation 
$$z(1-z)\frac{d^2w}{dz^2}  + \big[c-(a+b+1)z\big]\frac{dw}{dz}  - abw= 0.$$
Around the singular point $z = 0$, there are two independent solutions. One of them, if $c$ is not a non-positive integer, is
  \begin{equation}\label{hyper-2}
   w=w(a,b;c;z) = { }_{2}F_{1}(a,b;c;z) = \frac{1}{B(b,c-b)}\int_0^1x^{b-1}(1-x)^{c-b-1}(1-zx)^{-a}dx,
   \end{equation}
where $B(x,y)$ is the Beta function defined by
$$B(x,y)=\int_0^1t^{x-1}(1-t)^{y-1}dt$$
where  $x$ and $y$ are  complex numbers with positive  real part.

In fact, \eqref{hyper-2} converges for $z \in {\bf C}$ satisfying $\vert z\vert < 1$, and for $z = -1$ the definition is formally, however, it should be mentioned that for $z = -1$ and special $b$, $c$, and $a = t > 0$ the convergence of \eqref{hyper-2} makes sense.  

Several properties, as for example symmetries and some asymptotic behavior on the parameter $z$, appear on the literature about this function \eqref{hypergeomet}--\eqref{hyper-2} for particular cases of $a, b, c$. For example,  Bessel functions can be expressed as a limit of hypergeometric functions. However, there seems to be no any results on asymptotic behavior for the parameter $a$ when $a=t$ represents the time. On hypergeometric functions we can mention the works \cite{Erd}, \cite{Ga-ra}, \cite{Goursat}, \cite{Mill} and their  references. 

In this work we get some results for some fixed $a,b, c, z$. Our result seems new and very important. In our Theorem \ref{general-p}  we prove that 
$$I_p(t)  = \int_0^1(1+r^2)^{-t}r^pdr \sim t^{-\frac{p+1}{2}}, \quad t \gg 1$$
for each $p \geq 0$.
But we can note that by a change of variable
$$I_p(t)=\frac{1}{2}\int_0^1(1+s)^{-t}s^{\frac{p-1}{2}}ds 
= \frac{1}{2} B\Big(\frac{p+1}{2},1\Big){}_2F_1\Big(t, \frac{p+1}{2}; \frac{p+3}{2}; -1\Big),$$
since we choose in \eqref{hyper-2}: $b-1=\displaystyle{\frac{p-1}{2}}$, $c=b+1$, $a=t$ and $z=-1$. 
In this case,  one has
$${}_2F_1\Big(t, \frac{p+1}{2}; \frac{p+3}{2}; -1\Big)=(p+1)I_p(t),$$
because $ B\Big(\displaystyle{\frac{p+1}{2}},1\Big)=\displaystyle{\frac{2}{p+1}}$.

Then we obtain the following asymptotic behavior for a particular class of hypergeometric functions.
\begin{pro} \; Let $p\geq 0$. Then
\[{}_2F_1\Big(t, \frac{p+1}{2}; \frac{p+3}{2}; -1\Big) \sim t^{-\frac{p+1}{2}}, \quad t \gg 1.\]
\end{pro}

From this proposition we have in particular 
$${}_2F_1(t,\frac{1}{2};\frac{3}{2};-1)   \sim t^{-\frac{1}{2}}, \quad t \gg 1,$$
and 
$${}_2F_1(t,\frac{3}{2};\frac{5}{2};-1)  \sim t^{-\frac{3}{2}}, \quad t \gg 1.$$

As mentioned above, in the next subsection we show optimal asymptotic behavior of the hypergeometric functions given by the integral
$$I_p(t)= \int_0^1(1+r^2)^{-t}r^pdr, \quad t>1/2,$$
for each fixed  $p \geq 0$, and to the case for each $p \in {\bf R}$
$$J_p(t)=\int_1^{\infty}(1+r^2)^{-t}r^pdr, \quad t>1/2.$$
In particular,  it is known that
$$H_0(t):=\int_0^{\infty}(1+r^2)^{-t}dr=\frac{\sqrt{\pi}}{2}\frac{\Gamma(t-1/2)}{\Gamma(t)}, \quad t>1/2,$$
where $\Gamma=\Gamma(t)$ is the gamma function.

Then, by combining our decay estimates in the next section with $I_p(t)$ and $J_p(t)$ one can obtain the following asymptotic behavior of the function $\Gamma(t-1/2)/\Gamma(t)$
$$\frac{\Gamma(t-1/2)}{ \Gamma(t) }\sim t^{-1/2}, \quad t \gg1,$$ whose result does not seem to be well-known, although it is simple to see that  $\Gamma(t-1)/\Gamma(t)= \displaystyle{\frac{1}{t-1}}, t>1.$

\vspace{0.1cm}
Finally, it is important to observe that the behavior of hypergeometric functions of the type  ${\displaystyle \,_{2}F_{1}(t,b;c;-1)}$ appears when we study the asymptotic behavior of solutions for the wave equation under effects of  a special  dissipative term of logarithm type.

\subsection{General case}
Let $p \geq 0$ be a real number and  $I_p(t)$ be the function defined by 
 $$I_p(t)= \int_{0}^{1}(1+r^{2})^{-t}r^{p}dr, \quad \mbox{for } t>p.$$

The following theorem gives the optimal  asymptotic behavior of $I_p(t)$ for large $t$.
\begin{theo}\label{infit-1} Assume that $0\leq p \leq 3$. Then
$$I_p(t) \sim t^{-\frac{p+1}{2}}, \quad t \gg 1.$$
\end{theo}

{\it Proof.} Let $f(r)$ be the function given by $f(r)=(1+r^2)^{-t}r^p$, $\;r\geq 0$. Then $\beta=\sqrt{\displaystyle{\frac{p}{2t-p}}}$, $t>p$,  is a global maximum of $f$ and $0<\beta<1$ . Moreover $f(r)$ is a decreasing function for $r>\beta$,  increasing for $0<r<\beta$ when $p>0$,  $f(0)=1$ in case $p=0$ and $f(0)=0$ if $p>0$.

\vspace{0.2cm}
{\bf Case $p=0$:}
To prove this case  we  split the interval of  integration   in two  parts as follows.
 $$I_0(t)= \int_{0}^{t^{-1/2}}(1+r^{2})^{-t}dr
 + \int_{t^{-1/2}}^{1}(1+r^{2})^{-t}dr.$$

Now we note that 
\begin{align}\label{po1}
  \int_{0}^{t^{-1/2}}(1+r^{2})^{-t}dr  \leq t^{-1/2},
 \end{align}
because $f(0)=1$ is maximum global of $f(r)$ on the interval $(0,\infty)$.

On the other hand, by using a change of variable $u = \log(1+r^2)$ we have
\begin{align}\label{po2}
 \int_{t^{-1/2}}^{1}(1+r^{2})^{-t}dr&= \int_{t^{-1/2}}^{1}e^{-t\;\log(1+r^{2})}dr\nonumber\\
 &= \frac{1}{2} \int_{\log(1+\frac{1}{t})}^{\log2}e^{-(t-1)u}(e^u -1)^{-1/2}du\nonumber\\
&\leq   \frac{1}{2} \int_{\log(1+\frac{1}{t})}^{\log2}e^{-(t-1)u}(e^{\log(1+\frac{1}{t})} -1)^{-1/2}du\\
&\leq   \frac{1}{2} \int_{\log(1+\frac{1}{t})}^{\log2}e^{-(t-1)u}\;{t} ^{1/2}du\nonumber\\
&\leq  \frac{t ^{1/2}}{2}\;\frac{ e^{-(t-1)u}}{(t-1)}\Big|_{\log(1+\frac{1}{t})} \leq C\;\frac{t^{1/2}}{t-1}, \quad t\geq 2 \nonumber
\end{align}
with  $C>0$ a constant because  $ e^{-(t-1)\log(1+\frac{1}{t})}$ is a time-bounded function on $[2, \infty)$.

 \vspace{0.2cm}
The above estimates  give an optimal upper bound to  $I_0(t)$. 
 \vspace{0.2cm}
 
The  estimate to $I_0(t)$  from below is very easy.  Indeed, it  is obvious that  for $t>1$
$$
I_0(t) \geq \int_{0}^{t^{-1/2}}(1+r^{2})^{-t}dr \geq f(t^{-1/2})({t^{-1/2}} - 0)  = 
(1+ 1/t)^{-t}\;t^{-1/2}.
$$
Then, from the fact that  $\displaystyle{\lim_{t \rightarrow +\infty}}(1 + 1/t)^ {-t} = e^{-1}$, we  may fix arbitrary  positive  $C_0 < e^{-1}$ and choose $t_0>1$ depending on $C_0 $ such that 
\begin{align}\label{po3}
I_0(t) \geq C_0 \;t^{-1/2}, \quad t \geq t_0.
\end{align}

The estimates \eqref{po1}, \eqref{po2} and   \eqref{po3} prove the theorem to the case $p=0$. 
 \vspace{0.3cm}
 
{\bf Case $p>0$: }
To prove the theorem  for  $0< p \leq 3$  we  split the interval of  integration   in three parts, that is, we may write

\begin{align}\label{p1}
I_p(t)= \int_{0}^{\beta}(1+r^{2})^{-t}r^{p}dr + \int_{\beta}^{t^{-1/4}}(1+r^{2})^{-t}r^{p}dr 
 + \int_{t^{-1/4}}^{1}(1+r^{2})^{-t}r^{p}dr.
\end{align}
Note that $ \beta < t^{-1/4}$ for $t>p^2$.

\vspace{0.2cm}
The next step is to estimate each one of these integrals.
\vspace{0.2cm}
Based on the properties of $f(r)$  and the definition of $\beta$ we have 
\begin{align}\label{p11}
 \int_{0}^{\beta}(1+r^{2})^{-t}r^{p}dr &\leq f(\beta)(\beta-0)  = (1 + \beta^2)^{-t}\beta^{p+1}\\
 & =( 1 + \frac{p}{2t-p})^{-t} ( \frac{p}{2t-p})^{\frac{p+1}{2}}
 \leq C_p \Big(\frac{1}{2t-p}\Big)^{\frac{p+1}{2}}, \quad t>p, \nonumber
\end{align}
where $C_p>0 $ is a constant  depending on $p$ and we have used the fact that $( 1 + \displaystyle{\frac{p}{2t-p}})^{-t} $ is a time-bounded function  on the interval $[p, \infty)$.

Now we want to get an upper bound to the second integral on the right hand side of \eqref{p1}. To do that we perform the following estimates using the definition of $\beta$ and integration by parts.

\begin{align*}
\int_{\beta}^{t^{-1/4}}&(1+r^{2})^{-t}r^{p}dr = \frac{1}{2}\int_{\beta}^{t^{-1/4}}(1+r^{2})^{-t}2r r^{p-1}dr\\
&=   \frac{(1+r^2)^{-t+1}r^{p-1}}{2(-t+1)}\Big|_{\beta}^{t^{-1/4} } - 
\frac{1}{2}\int_{\beta}^{t^{-1/4}}\frac{(1+r^{2})^{-t+1}}{-t+1}(p-1)r^{p-2}dr\\
&\leq 
  \frac{(1+\frac{p}{2t-p})^{-t+1}(\frac{p}{2t-p})^{\frac{p-1}{2}}}{2(t-1)}
  + \frac{p-1}{2}\int_{\beta}^{t^{-1/4}}\frac{(1+r^{2})^{-t+1}}{t-1}r^{p-2}dr\\
  & \leq \frac{C_p}{(t-1){(2t-p)^{\frac{p-1}{2}}}} + \frac{p-1}{2(t-1)}\int_{\beta}^{t^{-1/4}}(1+r^{2})^{-t+1}r^{p-2}dr, \quad  t>\max\{1,p, p^2 \}. \\
\end{align*}

At this point we apply a change of variable $u=\log(1+ r^2)$ to obtain
\begin{align}\label{I2}
\int_{\beta}^{t^{-1/4}}&(1+r^{2})^{-t}r^{p}dr
  & \leq \frac{C_p}{(t-1){(2t-p)^{\frac{p-1}{2}}}} + \frac{p-1}{4(t-1)}\int_{\log(1+\beta^2)}^{\log(1+t^{-1/2})}e^{-(t-1)u}(e^u-1)^{\frac{p-3}{2}}du 
\end{align}
for all $t>max\{1,p\}.$

Now we also need to get an upper bound to  the integral on the right hand side of the above estimate for $0 \leq p \leq 3$.

\vspace{0.2cm}
For  $0 < p  \leq 3$ we may  estimate for $t> max\{1, p, p^2\}$

\begin{align}\label{I2pp}
  \frac{p-1}{4(t-1)}\int_{\log(1+\beta^2)}^{\log(1+t^{-1/2})}&e^{-(t-1)u}(e^u-1)^{\frac{p-3}{2}}du \nonumber\\
&\leq \frac{p-1}{4(t-1)}\int_{\log(1+\beta^2)}^{\log(1+t^{-1/2})}e^{-(t-1)u}(e^{\log(1+\beta^2)}-1)^{\frac{p-3}{2}}du \nonumber\\
&= \frac{(p-1)\beta^{p-3}}{4(t-1)}\int_{\log(1+\beta^2)}^{\log(1+t^{-1/2})}e^{-(t-1)u}du \\
&\leq \frac{(p-1)(\frac{p}{2t-p})^{\frac{p-3}{2}}}{4(t-1)}\;\frac{e^{-(t-1)u}}{t-1}\Big|_{\log(1+\beta^2)}\nonumber\\
&\leq \frac{C_p}{(t-1)^2(2t-p)^{\frac{p-3}{2}}} \leq \frac{C_p}{(t-1)^{\frac{p+1}{2}}} .\nonumber
\end{align}

The last above inequality with $C_p>0$  is due to the fact that the function $$e^{-(t-1)u}\Big|_{\log(1+\beta^2)}=(1+ \frac{p}{2t-p})^{-t+1}$$
is a bound function for $t>p$.
\vspace{0.3cm}

Next we need to estimate the third integral on the right hand side of \eqref{p1}. From the decreasing property of $f(r)$ on the interval of integration, we have 

\begin{align}\label{I3p}
 \int_{t^{-1/4}}^{1}(1+r^{2})^{-t}r^{p}dr \leq f(t^{- 1/4})(1-t^{-1/4}) \leq f(t^{- 1/4})=(1+\frac{1}{t^{1/2}})^{-t}t^{-\frac{p}{4}}.
\end{align}

Now we  observe that $\displaystyle{\lim_{t \rightarrow \infty}}(1+ \displaystyle{\frac{1}{\sqrt{t}}})^{-\sqrt{t}} = e^{-1}$. Then, there exists $t_0 >0$ such that
$$(1+ \frac{1}{\sqrt{t}})^{-\sqrt{t}} \leq  \frac{2}{e}, \quad t \geq t_0.$$
In particular 
\begin{align*}
(1+ \frac{1}{\sqrt{t}})^{-t} \leq \left(\frac{2}{e}\right)^{\sqrt{t}}\; = \left(\frac{e}{2}\right)^{-\sqrt{t}}\;, \quad t \geq t_0.
\end{align*}
Then combining this inequality with \eqref{I3p} we may conclude that 

\begin{align}\label{I3}
 \int_{t^{-1/4}}^{1}(1+r^{2})^{-t}r^{p}dr \leq \;t^{-\frac{p}{4}}\; \left(\frac{e}{2}\right)^{-\sqrt{t}}, \quad  t \geq t_0.
\end{align}

By substituting the estimates  \eqref{p11}, \eqref{I2} combined with  \eqref{I2pp} and \eqref{I3} in \eqref{p1},  we obtain the following optimal upper bound to $I_p(t)$ to the case $0<p \leq 3$.

\begin{align}\label{Ip-upper}
I_p(t) \leq C \; t^{-\frac{p+1}{2}}, \quad  t \gg 1. 
\end{align}

\vspace{0.3cm}

Finally we have to prove the lower estimate for the case $0<p \leq 3$.
To this case the function $f(r)=(1+r^2)^{-t}r^p$ is increasing on the interval $(0,\beta)$ with $\beta=\sqrt{\displaystyle{\frac{p}{2t-p}}}, \; t>p$. 

The  next estimate  give us the conclusion of  the proof of  the optimality of  the decay rate for $I_p(t), \; p>0$. 
In fact, the limit
\[\lim_{t \to \infty}(1 +  \frac{p}{4(2t-p)})^{-t}=e^{-\frac{p}{8}}\]
implies the existence of $t_0>0$ such that
$$I_p(t) \geq f\Big(\frac{\beta}{2}\Big)\Big(\beta-\frac{\beta}{2}\Big)
= \left(1+\frac{\beta^2}{4}\right)^{-t} \frac{\beta^{p+1}}{2^{p+1}} 
\geq \frac{1}{2}e^{-\frac{p}{8}} \frac{\beta^{p+1}}{2^{p+1}}
\geq \frac{e^{-\frac{p}{8}}}{2^{p+2}}\Big(\frac{p}{2t-p}\Big)^{\frac{p+1}{2}}, \quad t>p,
$$
because $\beta=\sqrt{\displaystyle{\frac{p}{2t-p}}}$ is the global  maximum of $f(r)$.
\vspace{0.2cm}

\hfill
$\Box$
\vspace{0.2cm}

The proof of Theorem 2.1 is now established. However,  our main aim in this section is to extend the result of this theorem for all $p \geq 0$. In order to do this we need the next important property of the  hypergeometric function 
$I_p(t)$ for $p \geq 2$.

\begin{lem}[Recurrence formula]\label{Rec-for}
Let $p \geq 2$ be a real number. Then 
$$I_p(t)=\frac{2^{-t+1}}{p+1-2t}  + \frac{p-1}{2t-p-1}I_{p-2}(t), \quad  t > \frac{p+1}{2}.$$
\end{lem}

{\it Proof.}\, Let $p \geq 2$, and $t > \displaystyle{\frac{p+1}{2}}$. It follows from integration by parts that
\begin{align*}
I_{p}(t) &= \frac{1}{2}  \int_{0}^{1}(1+r^{2})^{-t}2rr^{p-1}dr\\
&= \frac{2^{-t} }{1-t}+ \frac{p-1}{2t-2}\int_{0}^{1}(1+r^2)^{-t}(1+r^2)r^{p-2}dr   \\
&=  \frac{2^{-t} }{1-t}+ \frac{p-1}{2t-2}\int_{0}^{1}(1+r^2)^{-t}r^{p-2}dr  + \frac{p-1}{2t-2}\int_{0}^{1}(1+r^2)^{-t}r^{p}dr\\
 &=  \frac{2^{-t} }{1-t}+ \frac{p-1}{2t-2}I_{p-2}(t)  + \frac{p-1}{2t-2}I_p(t), 
 \end{align*}
which implies the identity
\[I_{p}(t) =  \frac{2^{-t+1} }{p+1- 2t}+  \frac{p-1}{2t-p-1} I_{p-2}(t).\] 

This yields the desired equality for $p \geq 2$. 
\hfill
$\Box$

Combining the recurrence formula with Theorem \ref{infit-1} we may prove the general result for $I_p(t)$.

\begin{theo}\label{general-p}
 Let $p \geq 0$ be a real number. Then
$$I_p(t) \sim t^{-\frac{p+1}{2}}, \quad t \gg 1.$$
\end{theo}

{\it Proof.}\,Applying Lemma \ref{Rec-for} for  $3\leq p \leq 4$ and using the result of Theorem \ref{infit-1}, which holds for $1\leq p-2 \leq 2$ we get the proof for $3\leq p \leq 4$. By a similar argument to the case for $4\leq p \leq 5$ and 
 $2\leq p-2 \leq 3$ we obtain the statement for $4\leq p \leq 5$. The general result follows  using the principle of induction.

\hfill
$\Box$

\begin{rem}{\rm It follows from Theorem \ref{general-p} that the optimal rate of decay of the function $I_{n-1}(t)$ is the same as that of the Gauss kernel in $L^{2}$-sense: $\Vert G(t,\cdot)\Vert^{2} \sim t^{-\frac{n}{2}}$ as $t \to \infty$, where 
\[G(t,x) := \frac{1}{(\sqrt{4\pi t})^{n}}e^{-\frac{\vert x\vert^{2}}{4t}}.\]}
\end{rem}

In order to deal with the high frequency part of estimates, one defines a function
$$J_p(t)=\int_1^{\infty}(1+r^2)^{-t}r^p dr$$
for $p \in {\bf R}$.

Then the next lemma is important to get estimates on the zone of high frequency to problem \eqref{eqn}--\eqref{initial}.
\begin{lem}\label{infit}
\,Let $p \in {\bf R}$. Then it holds that 
$$J_p(t) \sim \dfrac{2^{-t}}{t-1}, \quad t \gg 1.$$
\end{lem}
{\it Proof.}\,We first note that 
$$J_p(t)=\int_1^{\infty} e^{-t\log(1+r^2)}r^p dr, \quad  t>1.$$
Applying a change of variable $u=\log(1+r^2)$ we get 
$$J_p(t)=\frac{1}{2}\int_{\log2}^{\infty} e^{-(t-1)u}(e^u-1)^{\frac{p-1}{2}}du.$$
For $p < 1 $ and $u \geq \log 2$ we have
$$e^{\frac{p-1}{2}u}\leq (e^u-1)^{\frac{p-1}{2}} \leq\; 1.$$
Then using this inequality we obtain for $t>1$ the  double below-above estimate
$$
2^{\frac{p-1}{2}}\frac{2^{-t}}{t-\frac{p+1}{2}} = \frac{1}{2}\int_{\log2}^{\infty} e^{-(t-\frac{p+1}{2})u}du
\leq J_p(t)\leq\frac{1}{2}\int_{\log2}^{\infty} e^{-(t-1)u}du = \dfrac{2^{-t}}{t-1}
$$
For $p\geq 1 $ and $u \geq \log 2$ we have the inequality
$$1\leq (e^u-1)^{\frac{p-1}{2}} \leq e^{\frac{p-1}{2}u}.$$
Thus we also obtain for this case and $t>1$
$$ \dfrac{2^{-t}}{t-1} \leq J_p(t)\leq 2^{\frac{p-1}{2}}\frac{2^{-t}}{t-\frac{p+1}{2}}.$$
These estimates imply the lemma.
\hfill
$\Box$

\vspace{0.2cm}
For later use we prepare the following simple lemma, which implies the exponential decay estimates of the middle frequency part.
\begin{lem}\label{intermid}\,Let $p \in {\bf R}$, and $\eta \in (0,1]$. Then there is a constant $C > 0$ such that 
$$\int_{\eta}^{1}(1+r^{2})^{-t}r^{p}dr \leq C(1+\eta^{2})^{-t}, \quad t \geq 0.$$
\end{lem}


\subsection{ Inequalities and asymptotics }
\begin{lem}\label{ab-estim}
 Let $a(\xi) $ and $b(\xi) $ be the  functions given by
 \begin{equation}\label{roj}
a(\xi)=\dfrac{\log(1+|\xi|^{2})}{2} \quad \text{  and  }  \quad b(\xi)= \frac{1}{2}\sqrt{4|\xi|^{2}-\log^2(1+|\xi|^{2}) }
\end{equation}
for $\xi \in {\bf R}^{n}$. Then, the following estimates hold.\\
\vspace{0.2cm}
 $(i)\;\;\;\;\;\dfrac{|a(\xi)|^2}{|b(\xi)|^2} \leq \displaystyle{\frac{1}{3}}, \quad \xi \ne 0;$\\
\vspace{0.35cm}
$(ii)\;\;\; \dfrac{(b(\xi)-|\xi| )^2}{b(\xi)^2} \leq \displaystyle{\frac{28}{3}}, \quad \xi \ne 0.$
\end{lem}
{\it Proof.} To prove the lemma we use the elementary inequality
$$|\xi| - \log( 1+ |\xi|^2) \geq 0, \quad \mbox{for all}\; \xi \in {\bf R}^n.$$
Then 
\begin{align*} (ii) \;\;\dfrac{|a(\xi)|^2}{|b(\xi)|^2} =\dfrac{\log^2( 1+ |\xi|^2)}{4|\xi|^2 -\log^2( 1+ |\xi|^2)}
\leq \dfrac{ |\xi|^2}{4|\xi|^2 - |\xi|^2} \leq \; \displaystyle{\frac{1}{3}}, \quad \xi \ne 0 ;
\end{align*}
\begin{align*} (ii) \;\;\dfrac{(b(\xi)-|\xi| )^2}{b(\xi)^2}
&\leq 4(1 + \dfrac{ 4|\xi|^2}{b(\xi)^2})  \leq 4\left(1 + \dfrac{ 4|\xi|^2}{4|\xi|^2 -\log^2( 1+ |\xi|^2)}\right)\\
&\leq  4 (1+  \dfrac{4|\xi|^2}{3|\xi|^2}) \leq \dfrac{28}{3}, \quad \xi \ne 0.
\end{align*}
\hfill
$\Box$

In the next section, to study an asymptotic profile of the solution to problem \eqref{eqn}--\eqref{initial} we consider a decomposition of the Fourier transformed initial data.

\begin{rem}\label{obs1}
{\rm Using  the  Fourier transform we can get a decomposition of the initial data $\hat{u}_1$  as follows 
$$\hat{u}_1(\xi)=A_1(\xi)-iB_1(\xi)+P_1,\quad \xi \in {\bf R}^n, $$
where  $P_1, A_1, B_1$  are defined by  
\[P_{1} = \int_{{\bf R}^n}u_{1}(x) dx, \quad A_{1} (\xi)=\int_{{\bf R}^n}u_{1}(x)\big(1-\cos(\xi x) \big)dx, \quad B_{1}(\xi) =\int_{{\bf R}^n}u_{1}(x)\sin(\xi x)dx.\]}
\end{rem}

According to the above decomposition we can derive the following lemma (see Ikehata \cite{I-04}).
 
\begin{lem}\label{lema2.5.4}

Let $\kappa \in [0,1]$. For $u_{1} \in L^{1,\kappa}({\bf R}^{n})$ and $\xi \in {\bf R}^{n}$ it holds that 
$$|A_1(\xi)|\leq K|\xi|^\kappa\|u_{1}\|_{L^{1,\kappa}} \quad \text{ and } \quad |B_1(\xi)|\leq M|\xi|^\kappa\|u_{1}\|_{L^{1,\kappa}},$$
with positive constants $K$ and $M$ depending only on $n$.\\
\end{lem}


In order to show the optimality of the decay rates we need next two lemmas.
\begin{lem}\label{2.8}
Let  $n> 2$. Then there exists $t_0 > 0$ such that  for $t \geq t_{0}$ it holds that
$$C^{-1} t^{-\frac{n-2}{2}} \geq \int_{{\bf R}^{n}}{e^{-t \log(1+ |\xi|^2) }\dfrac{|\sin(|\xi|t)|^2}{|\xi|^2}}d\xi \geq C t^{-\frac{n-2}{2}}, $$
with  $C$ a positive constant depending only on   $n$.
\end{lem}
{\it Proof.}\,  First, we may note that 
\begin{align*}
M(t): & = \int_{{\bf R}^{n}} e^{-t\; \log(1+|\xi|^2)} \dfrac{|\sin(|\xi|t)|^2}{|\xi|^2} d\xi \\
&=\omega_n \int_0^{\infty}e^{-t\; \log(1+r^2)}  r^{n-3}|\sin(rt)|^2   dr\\
& \geq\omega_n \int_0^{\infty}e^{-t\;r^2}  r^{n-3}|\sin(rt)|^2   dr.\\
\end{align*}
Considering a change of variable $s=r\sqrt{t}$, for fix $t>0$  we arrive at
\begin{align*}
M(t)& \geq \dfrac{ \omega_n}{t^{\frac{n-2}{2}}}\int_0^{\infty}e^{-s^{2}}s^{n-3}\sin^2\big(s\sqrt{t})ds.
\end{align*}
\noindent
Using the identity 
\[2\sin^2 x = \left(1-\cos 2x\right),\]
we obtain
\begin{align*}
M(t) &\geq \frac{1}{2}\omega_n t^{-\frac{n-2}{2}} \int_0^{\infty}e^{-s^{2}}s^{n-3}\Big(1-\cos(2s\sqrt{t})\Big)ds\\
&=\frac{1}{2} \omega_n t^{-\frac{n-2}{2}} \big(A-F_n(t)\big),
\end{align*}
where
\[A = \displaystyle \int_0^{\infty}e^{-s^{2}}s^{n-3}ds, \quad \displaystyle  F_n(t)=\int_0^{\infty}e^{-s^{2}}s^{n-3}\cos\big(2s\sqrt{t}\big)ds.\]
Due to  the fact $e^{-s^{2}}s^{n-3}\; \in L^1({\bf R})$\; ($n>2$),  we can apply the Riemann-Lebesgue theorem to get 
$$F_n(t) \to 0, \quad t \to \infty .  $$
\noindent
Then we conclude the existence of  $t_0 >0$  such that  $F_n(t) \leq \displaystyle{\frac{A}{2}}$ for all $t \geq t_0$. Thus, the half part of lemma is proved with $C=\dfrac{\omega_nA}{4}$.

Next, let us prove upper bound of decay estimates. Indeed, 
\begin{align*}
M(t) &\leq \int_{{\bf R}^{n}} e^{-t\; \log(1+|\xi|^2)}\vert\xi\vert^{-2} d\xi \\
&=\omega_n \int_0^{\infty}e^{-t\; \log(1+r^2)}  r^{n-3} dr\\
&=\omega_{n}\int_{0}^{\infty}(1+r^{2})^{-t}r^{n-3}dr = \omega_{n}I_{n-3}(t) + \omega_{n}J_{n-3}(t)\\
&\leq C_{1,n}t^{-\frac{n-2}{2}} + C_{2,n}\frac{2^{-t}}{t},
\end{align*}
where one has just used Theorem  \ref{general-p}  and Lemma \ref{infit}. 
These imply the desired estimates.
\hfill
$\Box$

Following the same ideas of Lemma \ref{2.8} one can prove the following result, however, this is not used in the paper.
\begin{lem}\label{2.9}
Let  $n\geq 1$ . Then there exists  $t_0> 0$  such that   vale 
$$C_n^{-1} t^{-\frac{n}{2}} \geq \int_{{\bf R}^{n}}e^{-t\;\log(1+|\xi|^{2})}\cos^2(|\xi|t) d\xi \geq C_n t^{-\frac{n}{2}},  \quad t\geq t_0, $$
where $C_n$ is a positive constant depending only on  $n$.
\end{lem}


\section{Asymptotic profiles of solutions}

The  associated Cauchy problem to \eqref{eqn}-\eqref{initial} in the  Fourier space is given by 
\begin{align}\label{uhat}
&\hat{u}_{tt}(t,\xi) + |\xi|^2\hat{u}(t,\xi) + \log(1+  |\xi|^2)\hat{u}=0,\\
&\hat{u}(0,\xi)=u_0(\xi), \quad
\hat{u_t}(0,\xi)=u_1(\xi)\nonumber.
\end{align}

The characteristics roots  $\lambda_+$ and  $\lambda_-$  of the  characteristic polynomial 
$$\lambda^2+  \log(1+|\xi|^2)\lambda +|\xi|^{2}=0, \quad  \xi \in {\bf R}^n$$
associated to the equation \eqref{uhat}
 are given by 
\begin{align}\label{7.1}
\lambda_{\pm} = \dfrac{-\log(1+|\xi|^{2}) \pm \sqrt{\log^2(1+|\xi|^{2}) - 4| \xi|^{2}}}{2}.
\end{align}
It should be mentioned that  $\log(1+|\xi|^2) -  4| \xi|^{2} <0 $  for all $\xi \in {\bf R}^n, \xi \ne 0$, and the characteristics roots are complex  and the real part is negative, for all  $\xi \in {\bf R}^n, \xi \ne 0$. Then we can write down $\lambda_{\pm}$ in the following form
$$\lambda_\pm = - a(\xi)\pm i b(\xi),$$
where $a(\xi)$ and $b(\xi)$ are defined by \eqref{roj} in Lemma \ref{ab-estim}. In this case the solution of the problem \eqref{uhat} is given explicitly by
$$\hat{u}(t,\xi)=\left(\hat{u}_{0}(\xi)\cos(b(\xi)t) + \dfrac{\hat{u}_1(\xi) +\hat{u}_0(\xi) a(\xi)}{b(\xi)}\sin(b(\xi)t)\right)e^{-a(\xi)t}$$
for  $\xi \in {\bf R}^n, \xi \ne 0 $ and $t \geq 0$.

Next, in order to find a better expression for $\hat{u}(t,\xi)$ we apply the mean value theorem to get

\begin{equation}\label{sine}
\sin\left(b(\xi)t\right)=\sin(|\xi|t)+t\left(b(\xi)-|\xi|\right)\cos(\mu(\xi)t),
\end{equation}
with
$$\mu(\xi):=\theta_{1} b(\xi)+(1-\theta_{1})|\xi|$$
for some $\theta_{1} \in (0,1)$, and

\begin{equation}\label{root}
\frac{1}{\sqrt{1-g(r)}} = 1+\frac{\log^{2}(1+r^{2})}{8r^{2}}\frac{1}{\sqrt{(1-\theta_{2}g(r))^{3}}}
\end{equation}
with some $\theta_{2} \in (0,1)$, where $r := \vert\xi\vert$, and
\[g(r) := \frac{\log^{2}(1+r^{2})}{4r^{2}}.\]
\noindent
The identity \eqref{root} was obtained applying the mean value theorem to the function $$G(s)=\frac{1}{\sqrt{(1-sg(r))^{3}}}, \; 0 \leq s \leq 1.$$
\noindent
Then by using Remark \ref{obs1}, \eqref{sine} and \eqref{root} $\hat{u}(t,\xi)$ can be re-written as 

\begin{align}\label{express1}
\hat{u}(t,\xi) &= P_{1}e^{-a(\xi)t}\frac{\sin(tr)}{r} + P_{1}\frac{\log^{2}(1+r^{2})}{8r^{3}}\frac{1}{\sqrt{(1-\theta_{2}g(r))^{3}}}e^{-a(\xi)t}\sin(tr)\nonumber\\
&+e^{-a(\xi)t}\cos(b(\xi)t)\hat{u}_0(\xi) + \left(\frac{a(\xi)}{b(\xi)}\right)e^{-a(\xi)t}\sin(b(\xi)t)\hat{u}_0(\xi)\\
&+ \left(\frac{A_{1}(\xi) - iB_{1}(\xi)}{b(\xi)}\right)e^{-a(\xi)t}\sin(b(\xi)t) + P_{1}t e^{-a(\xi)t}\left(\frac{b(\xi)-r}{b(\xi)}\right)\cos(\mu(\xi)t).\nonumber
\end{align}

We want to introduce an asymptotic profile as $t \to \infty$ in a simple form:
\begin{equation}\label{leading}
P_{1}e^{-a(\xi)t}\dfrac{\sin(|\xi|t)}{|\xi|},
\end{equation}
where $a(\xi)=\displaystyle{\frac{\log(1+|\xi|^2)}{2}}$.

Our goal in this section is to get decay estimates in time to the remainder therms defined in \eqref{express1}. To proceed with that we define the next $5$ functions which imply remainders with respect to the leading term \eqref{leading}.
 \begin{itemize}
\item[$\bullet$] $K_1(t, \xi)=\Big(\displaystyle{\frac{A_1(\xi)-iB_1(\xi)}{b(\xi)}}\Big)e^{-a(\xi)t}\sin(b(\xi)t)$;
\item[$\bullet$] $K_2(t, \xi)= \hat{u}_{0}(\xi) \displaystyle{\frac{a(\xi)}{b(\xi)}}e^{-a(\xi)t}\sin\big( b(\xi)t\big)$;
\item[$\bullet$] $K_3(t, \xi)= \hat{u}_{0}(\xi) e^{-a(\xi)t}\cos\big( b(\xi)t\big)$;
\item[$\bullet$] $K_4(t, \xi)= P_{1}e^{-a(\xi)t}\sin(rt)\displaystyle{\frac{\log^{2}(1+r^{2})}{8r^{3}}}\displaystyle{\frac{1}{\sqrt{(1-\theta_{2}g(r))^{3}}}} ,  \quad r=|\xi|>0$;
\item[$\bullet$] $K_5(t, \xi)= P_{1}e^{-a(\xi)t}t\left(\dfrac{b(\xi)-|\xi|}{b(\xi)}\right)\cos(\mu(\xi)t)$,
 \end{itemize}
where $a(\xi)$ and $b(\xi)$ are defined in Lemma \ref{ab-estim}. Note that using these $K_{j}(t,\xi)$ ($j = 1,2,3,4,5$) the solution $\hat{u}(t,\xi)$ to problem \eqref{uhat} can be expressed as
\begin{equation}\label{expression}
\hat{u}(t,\xi) - P_{1}e^{-a(\xi)t}\frac{\sin(tr)}{r} = \sum_{j = 1}^{5}K_{j}(t,\xi).
\end{equation}

Let us check, in fact, that $\{K_{j}(t,\xi)\}$ become error terms by using previous lemmas studied in Section 2. 

First we obtain decay rates for each one of these functions on the zone of low frequency $|\xi| \ll 1$. 

We begin with $K_1(t,\xi)$.
\vspace{0.2cm}

For this function we prepare the following expression for $1/b(\xi)$ based on \eqref{root}:
\begin{equation}\label{rootb}
\frac{1}{b(\xi)} = \frac{1}{r} + \frac{\log^{2}(1+r^{2})}{8r^{3}}\frac{1}{\sqrt{(1-\theta_{2}g(r))^{3}}} , \quad r=|\xi| >0.
\end{equation}
Then, 
\begin{align*}
K_{1}(t,\xi) &: \frac{A_1(\xi)-iB_1(\xi)}{\vert\xi\vert}e^{-a(\xi)t}\sin(b(\xi)t) \\
&+ \big(A_1(\xi)-iB_1(\xi)\big)\frac{\log^{2}(1+r^{2})}{8r^{3}}\frac{e^{-a(\xi)t}\sin(b(\xi)t)}{\sqrt{(1-\theta_{2}g(r))^{3}}} =: K_{1,1}(t,\xi) + K_{1,2}(t,\xi).
\end{align*}
\noindent
It is easy to check the following estimate based on Lemma \ref{lema2.5.4} with $k = 1$ and 
Theorem  \ref{general-p}:
\begin{align}\label{K{1,1}}
\int_{\vert\xi\vert\leq 1}\vert K_{1,1}(t,\xi)\vert^{2}d\xi 
&\leq (M+K)^{2}\Vert u_{1}\Vert_{1,1}^{2}\int_{\vert\xi\vert\leq 1}e^{-t\log(1+\vert\xi\vert^{2})}d\xi  \nonumber\\
&\leq \omega_n (M+K)^{2}\Vert u_{1}\Vert_{1,1}^{2}\int_0^1(1+r^{2})^{-t}r^{n-1}dr \nonumber \\
&\leq C\omega_{n}(M+K)^{2}t^{-\frac{n}{2}}\Vert u_{1}\Vert_{1,1}^{2}, \quad (t \gg 1).
\end{align}

On the other hand, since
\[\lim_{r \to +0}\frac{\log^{2}(1+r^{2})}{r^{2}} = 0,\]
there is a constant $\delta > 0$ such that for all $0 < r \leq \delta$ it holds that
\begin{equation}\label{bound}
\frac{\log^{2}(1+r^{2})}{4r^{2}} \leq \frac{1}{2}.
\end{equation}
Then, the definition of $g(r)$ implies 
\begin{equation}\label{defforf}
\frac{1}{\sqrt{(1-\theta_{2}g(r))^{3}}} \leq 2\sqrt{2}.
\end{equation}
Thus,  from \eqref{bound}, \eqref{defforf} and Theorem  \ref{general-p}  together with Lemma \ref{lema2.5.4} for  $k = 1$, one has
\begin{align}\label{K{1,2}}
\int_{\vert\xi\vert\leq \delta}\vert K_{1,2}(t,\xi)\vert^{2}d\xi& \leq 8^{-1}(M+K)^{2}\Vert u_{1}\Vert_{1,1}^{2}\int_{\vert\xi\vert\leq \delta}\left(\frac{\log^{2}(1+r^{2})}{r^{2}}\right)^{2}e^{-2ta(\xi)}d\xi \nonumber\\
&\leq \frac{1}{2}(M+K)^{2}\Vert u_{1}\Vert_{1,1}^{2}\int_{\vert\xi\vert\leq \delta}e^{-2ta(\xi)}d\xi \nonumber\\
&\leq \frac{1}{2}(M+K)^{2}\Vert u_{1}\Vert_{1,1}^{2}\omega_{n}\int_{0}^{1}(1+r^{2})^{-t}r^{n-1}d\xi \nonumber\\
&\leq C\omega_{n}(M+K)^{2}\Vert u_{1}\Vert_{1,1}^{2}t^{-\frac{n}{2}}, \quad (t \gg 1).
\end{align}
By combining \eqref{K{1,1}} and \eqref{K{1,2}}  we have the  following estimate for $K_{1}(t,\xi)$,
\begin{equation}\label{K{1}}
\int_{\vert\xi\vert\leq\delta}\vert K_{1}(t,\xi)\vert^{2}d\xi \leq C_{1,n}\Vert u_{1}\Vert_{1,1}^{2}t^{-\frac{n}{2}}, \quad (t \gg 1).
\end{equation}
Similarly to the computation for \eqref{K{1}}, one can obtain the estimate for $K_{4}(t,\xi)$
\begin{equation}\label{K{4}}
\int_{\vert\xi\vert\leq\delta}\vert K_{4}(t,\xi)\vert^{2}d\xi \leq C_{1,n}\vert P_{1}\vert t^{-\frac{n}{2}}, \quad (t \gg 1),
\end{equation}
because 
\begin{equation}\label{flim}
\lim_{r \to +0}\frac{\log^{2}(1+r^{2})}{r^{3}} = 0.
\end{equation}
For $K_{2}(t,\xi)$ and $K_{3}(t,\xi)$, by using (i) of Lemma \ref{ab-estim} one can easily obtain the estimate:
\begin{equation}\label{K{j}}
\int_{\vert\xi\vert\leq 1}\vert K_{j}(t,\xi)\vert^{2}d\xi \leq C_{1,n}\Vert u_{0}\Vert_{1} t^{-\frac{n}{2}}, \quad (t \gg 1),
\end{equation} 
for each $j = 2,3$. So, it suffices to deal with the case for $K_{5}(t,\xi)$. For this we remark that
\[b(\xi) - r = r\left(-\frac{\frac{\log^{2}(1+r^{2})}{4r^{2}}}{1+\sqrt{1-\frac{\log^{2}(1+r^{2})}{4r^{2}}}}\right).\] 
This implies
\[\vert b(\xi) - r\vert \leq r^{3}\frac{\log^{2}(1+r^{2})}{4r^{4}} =: r^{3}h(r),\]
where we see
\[h(r) = \frac{\log^{2}(1+r^{2})}{4r^{4}} \to \frac{1}{4}\quad (r \to +0).\]
So, there exists a constant $\delta_{0} > 0$ such that for all $r \in (0,\delta_{0}]$ it holds that
\[0 < \frac{\log^{2}(1+r^{2})}{4r^{4}} \leq 1.\]
Thus, one can estimate $K_{5}(t,\xi)$ as follows:
\begin{equation}\label{delta}
\int_{\vert\xi\vert\leq \delta_{0}}\vert K_{5}(t,\xi)\vert^{2}d\xi \leq \vert P_{1}\vert^{2}t^{2}\int_{\vert\xi\vert\leq \delta_{0}}r^{6}\frac{e^{-2ta(\xi)}}{b(\xi)^{2}}d\xi.
\end{equation}
\noindent
On the other hand, since $\log^2(1+r^2) \leq 2r^2$ for all $r\geq 0$,  it follows that 
$$\frac{1}{b(\xi)^2} = \frac{4}{4r^2 - \log^2(1+r^2)} \leq \frac{2}{r^2}, \quad |\xi|= r>0.$$ 
Therefore, one can estimate \eqref{delta} for $r \in (0,\delta_{1}]$ with sufficiently small $\delta_{1} \leq \delta_{0}$ as follows
\begin{equation}\label{delta{1}}
\int_{\vert\xi\vert\leq \delta_{1}}\vert K_{5}(t,\xi)\vert^{2}d\xi \leq 2\vert P_{1}\vert^{2}t^{2}\int_{\vert\xi\vert\leq \delta_{1}}r^{4}\frac{e^{-2ta(\xi)}}{b(\xi)^{2}}d\xi \leq C\vert P_{1}\vert^{2}t^{-\frac{n}{2}},
\end{equation}
where one has just used Theorem \ref{general-p} and the definition of $a(\xi)$ in Lemma  \ref{ab-estim}.

\vspace{0.2cm}
Now, by summarizing above discussion one can arrived at the following crucial lemma based on \eqref{expression}, \eqref{K{1}}, \eqref{K{4}}, \eqref{K{j}}, and \eqref{delta{1}}.

\begin{pro}\label{proposition3.1}\, Let $n \geq 1$. Then, there exists a small constant $\delta_{1} \in (0,1]$ such that
\[\int_{\vert\xi\vert\leq \delta_{1}} \big|\hat{u}(t,\xi) -  P_{1}e^{-a(\xi)t}\frac{\sin(tr)}{r}\big|^{2}d\xi 
\leq C\big(\vert P_{1}\vert^{2} + \Vert u_{0}\Vert_{1}^{2} + \Vert u_{1}\Vert_{1,1}^{2}\big)t^{-\frac{n}{2}},\quad (t \gg 1),\]
with some generous constant $C=C_n > 0$ depending only on  the dimension $n$.  
\end{pro}

Next, let us prepare the so-called high frequency estimates for such error terms $K_i(t,\xi)$. These terms decay very fast, as usual.

\begin{lem}\label{lemma3.1}
Let $n \geq 1$. Then, it holds that
\[\int_{|\xi|\geq \delta_{1}}|K_i(t,\xi)|^2d\xi \leq C ||u_j||_{1}^{2} o(t^{-\frac{n}{2}}), \quad(t \to \infty),\]
where $j=1$ for $i=1, 4, 5$ and $j=0$ for $i=2,3$, and $\delta_{1} > 0$ is a number defined in Proposition {\rm \ref{proposition3.1}}.
\end{lem}
{\it Proof.}\, We give the proof only for $K_5(t,\xi)$. The other cases are similar. Indeed, it follows from (ii) of Lemma \ref{ab-estim}, Lemmas \ref{infit} and \ref{intermid} that
\begin{align*}
\int_{|\xi|\geq \delta_{1}}|K_5(t,\xi)|^2d\xi &\leq \frac{16}{3}\vert P_{1}\vert^{2} t^{2}\int_{|\xi|\geq \delta_{1}}e^{-\log(1+\xi^2)t}d\xi \\
&\leq \frac{16}{3}\vert P_{1}\vert^{2} t^{2}\left(\int_{1\geq |\xi|\geq\delta_{1}}e^{-\log(1+\xi^2)t}d\xi + \int_{|\xi|\geq\delta_{1}}e^{-\log(1+\xi^2)t}d\xi\right)\\
&\leq \frac{16}{3}\vert P_{1}\vert^{2} t^{2}\omega_{n}\left(\int_{\delta_{1}}^{1}(1+r^{2})^{-t}r^{n-1}dr + \int_{1}^{\infty}(1+r^{2})^{-t}r^{n-1}dr \right)\\
&\leq \frac{16}{3}\vert P_{1}\vert^{2} t^{2}\omega_{n}\left(\int_{\delta_{1}}^{1}(1+r^{2})^{-t}r^{n-1}dr + \int_{1}^{\infty}(1+r^{2})^{-t}r^{n-1}dr \right)\\
&\leq \frac{16}{3}\vert P_{1}\vert^{2} t^{2}\omega_{n}\left( C(1+\delta_{1}^{2})^{-t} + \frac{2^{-t}}{t-1}  \right),\end{align*}
which implies the desired estimate for $K_{5}(t,\xi)$.
\hfill
$\Box$

Now, as a direct consequence of Lemma \ref{lemma3.1} and \eqref{expression} one can get the high frequency estimates for the error terms.
\begin{pro}\label{proposition3.2}\, Let $n \geq 1$. Then, there exists a small constant $\delta_{1} > 0$ such that
\[\int_{\vert\xi\vert \geq \delta_{1}}\vert\hat{u}(t,\xi) - P_{1}e^{-a(\xi)t}\frac{\sin(t\vert\xi\vert)}{\vert\xi\vert}\vert^{2}d\xi \leq C(\Vert u_{1}\Vert_{1}^{2} + \Vert u_{0}\Vert_{1}^{2})o(t^{-\frac{n}{2}}),\quad (t \to \infty),\]
with some generous constant $C > 0$.   
\end{pro}
 
Finally, Theorem \ref{main-theo} is a direct consequence of Propositions \ref{proposition3.1} and \ref{proposition3.2}

\begin{rem}\,{\rm The decay rate stated in Proposition \ref{proposition3.2} can be drawn with a more precise fast decay rate, however, since the decay rate in Proposition \ref{proposition3.1} is essential, and the rate of decay in Proposition \ref{proposition3.2} can be absorbed into that of Proposition \ref{proposition3.1}, we have employed such style for simplicity.}
\end{rem} 


\section{Optimal rate of decay of solutions.}

In this section we study the optimal decay rate in the sense of $L^{2}$-norm of the solutions to problem (1.1)-(1.2).\\
We first prepare the following proposition in the one dimensional case. 

\begin{pro}\label{proposition4.1}\,It is true that
\[\int_{{\bf R}}(1+\vert\xi\vert^{2})^{-t}\frac{\sin^{2}(t\vert\xi\vert)}{\vert\xi\vert^{2}}d\xi \sim t,\quad (t \gg 1).\]
\end{pro}

{\it Proof.}\,We set
\[Q(t) := \int_{0}^{\infty}(1+r^{2})^{-t}\frac{\sin^{2}(tr)}{r^{2}}dr,\]
and it suffices to obtain the estimate stated in Proposition for $Q(t)$. Then, $Q(t)$ can be divided into two parts:
\[Q(t) = Q_{l}(t) + Q_{h}(t),\]
with
\[Q_{l}(t) := \int_{0}^{1/t}(1+r^{2})^{-t}\frac{\sin^{2}(tr)}{r^{2}}dr,\]
\[Q_{h}(t) := \int_{1/t}^{\infty}(1+r^{2})^{-t}\frac{\sin^{2}(tr)}{r^{2}}dr.\]
\noindent
{\bf (i)}\,upper bound for $Q_{j}(t)$ with $j = l,h$. \\
Indeed, let $1/t < 1$. Then, 
\[Q_{l}(t) \leq \int_{0}^{1/t}(1+r^{2})^{-t}\frac{(tr)^{2}}{r^{2}}dr = t^{2}\int_{0}^{1/t}(1+r^{2})^{-t}dr \leq t^{2}\int_{0}^{1/t}dr = t,\]
where one has just used the fact that if $0 \leq tr \leq 1$  then $0 \leq \sin(tr) \leq tr$. This implies
\begin{equation}\label{41}
Q_{l}(t) \leq Ct\quad (t > 1).
\end{equation} 
On the other hand, it follows from integration by parts one can get
\[Q_{h}(t) \leq \left[-r^{-1}(1+r^{2})^{-t}\right]_{r:=1/t}^{r:=\infty} - 2t\int_{1/t}^{\infty}(1+r^{2})^{-t-1}dr \leq t(1+\frac{1}{t^{2}})^{-t}.\]
Now, since
\begin{equation}\label{42}
\lim_{t \to \infty}(1+\frac{1}{t^{2}})^{-t} = 1,
\end{equation}
there is a constant $t_{0} \gg 1$ such that for all $t \geq t_{0}$
\begin{equation}\label{43}
Q_{h}(t) \leq 2t.
\end{equation}
\noindent
Estimates \eqref{41} and \eqref{43} imply
\begin{equation}\label{44}
Q(t) \leq Ct\quad (t \gg 1).
\end{equation}
\noindent
{\bf (ii)}\,lower bound for $Q_{j}(t)$ with $j = l,h$. \\
Indeed, let $1/t < 1$ again. Then, since $2\sin(tr) \geq tr$ if $0 \leq tr \leq 1$,  the following estimate holds for $t>1$
\[Q_{l}(t) \geq \frac{t^{2}}{4}\int_{0}^{1/t}(1+r^{2})^{-t}dr \geq \frac{t^{2}}{4}(1+\frac{1}{t^{2}})^{-t}\frac{1}{t} = \frac{t}{4}(1+\frac{1}{t^{2}})^{-t}.\]
Thus, because of \eqref{42}  we can get
\begin{equation}\label{45}
Q_{l}(t) \geq Ct\quad (t \gg 1).
\end{equation} 
\noindent
To treat $Q_{h}(t)$ we set
\begin{equation}\label{nyu}
\nu := \frac{5\pi}{4t},\qquad \nu' := \frac{7\pi}{4t},
\end{equation}
and 
\[\rho := \frac{49\pi^{2}}{16}.\]
If $\nu \leq r \leq \nu'$, then one has
\begin{equation}\label{sineestimate}
\vert\sin(tr)\vert \geq \frac{1}{\sqrt{2}}.
\end{equation}
So, one can get a series of estimates from below because of $\displaystyle{\frac{1}{t}} < \displaystyle{\frac{5\pi}{4t}}$ ($1 < t$):
\begin{align*}
Q_{h}(t) & \geq \frac{1}{2}\int_{\nu}^{\nu'}(1+r^{2})^{-t}r^{-2}dr 
\geq \frac{1}{2}(\frac{7\pi}{4t})^{-2}\int_{\nu}^{\nu'}(1+r^{2})^{-t}dr,\\  
&\geq \frac{1}{2}\frac{16t^{2}}{49\pi^{2}}(1+\frac{\rho^{2}}{t^{2}})^{-t}(\nu'-\nu)\\
&= \frac{4}{49\pi}(1+\frac{\rho}{t^{2}})^{-t}t.
\end{align*}
Since
\[\lim_{t \to \infty}(1+\frac{\rho}{t^{2}})^{-t} = 1,\]
one can arrive at the crucial estimate:
\begin{equation}\label{46}
Q_{h}(t) \geq Ct \quad (t \gg 1).
\end{equation} 

By combining \eqref{45} and \eqref{46} it results that
\begin{equation}
Q(t) \geq Ct \quad (t \gg 1).
\end{equation} 

Finally, the desired estimate can be accomplished by \eqref{44} and (4.9).
\hfill
$\Box$

Next we deal with the two dimensional case, which is rather difficult. 
\begin{pro}\label{proposition4.2}\,It is true that
\[\int_{{\bf R^{2}}}(1+\vert\xi\vert^{2})^{-t}\frac{\sin^{2}(t\vert\xi\vert)}{\vert\xi\vert^{2}}d\xi \sim \log t,\quad (t \gg 1).\]
\end{pro}
{\it Proof.}\,\,It suffices to get the result to the following function after polar coordinate transform:
\[R(t) := \int_{0}^{\infty}(1+r^{2})^{-t}\frac{\sin^{2}(tr)}{r}dr.\]
Then, $R(t)$ can be divided into two parts:
\[R(t) = R_{l}(t) + R_{h}(t),\]
with
\[R_{l}(t) := \int_{0}^{1/t}(1+r^{2})^{-t}\frac{\sin^{2}(tr)}{r}dr,\]
\[R_{h}(t) := \int_{1/t}^{\infty}(1+r^{2})^{-t}\frac{\sin^{2}(tr)}{r}dr.\]
\noindent
{\bf (i)}\,upper bound for $R_{j}(t)$ with $j = l,h$. 
\vspace{0.15cm}

Indeed, let $1/t < 1$. Then, 
\begin{align*}
R_{l}(t)& \leq \int_{0}^{1/t}(1+r^{2})^{-t}\frac{(tr)^{2}}{r}dr = t^{2}\int_{0}^{1/t}(1+r^{2})^{-t}r dr\\
&\leq t\int_{0}^{1/t}(1+r^{2})^{-t}dr \leq t \int_{0}^{1/t}dr = 1,
\end{align*}
where one has just used the fact that if $0 \leq tr \leq 1$  then $0 \leq \sin(tr) \leq tr$. This implies
\begin{equation}\label{47}
R_{l}(t) \leq 1\quad (t > 1).
\end{equation} 

On the other hand, it follows from the integration by parts one can get
\begin{align}\label{48}
R_{h}(t)& \leq \int_{1/t}^{\infty}(1+r^{2})^{-t}r^{-1}dr\nonumber\\
&= \left[(\log r)(1+r^{2})^{-t}\right]_{r:=1/t}^{r:=\infty} + 2t\int_{1/t}^{\infty}(r\log r)(1+r^{2})^{-t-1}dr\nonumber\\
&= (\log t)(1+ t^{-2})^{-t} + 2t\int_{1/t}^{1}(r\log r)(1+r^{2})^{-t-1}dr + 2t\int_{1}^{\infty}(r\log r)(1+r^{2})^{-t-1}dr\nonumber\\
&\leq (\log t)(1+ t^{-2})^{-t} + 2t\int_{1}^{\infty}(r\log r)(1+r^{2})^{-t-1}dr\nonumber\\
&\leq (\log t)(1+ t^{-2})^{-t} + 2t\int_{1}^{\infty}r^{2}(1+r^{2})^{-t-1}dr\\
&\leq (\log t)(1+ t^{-2})^{-t} + 2t\int_{1}^{\infty}(1+r^{2})(1+r^{2})^{-t-1}dr\nonumber\\
&
= (\log t)(1+ t^{-2})^{-t} + 2t\int_{1}^{\infty}(1+r^{2})^{-t}dr.\nonumber
\end{align} 

Now, because of Lemma 2.2 one can see that
\begin{equation}\label{49}
\int_{1}^{\infty}(1+r^{2})^{-t}dr \leq C\frac{2^{-t}}{t-1},\quad (t \gg 1),
\end{equation}
and since $\displaystyle{\lim_{t \to \infty}}(1+t^{-2})^{-t} =1$, one has $(1+t^{-2})^{-t} \leq 2$ for $t \gg 1$. Thus, it follows from \eqref{48} and \eqref{49}, one can arrive at the estimate:
\begin{equation}\label{410}
R_{h}(t) \leq 2\log t + C 2^{-t}, \quad (t \gg 1).
\end{equation}  
\eqref{47} and \eqref{410} implies the upper bound for $t$ of the quantity $R(t)$:
\begin{equation}\label{411}
R(t) \leq C\log t, \quad (t \gg 1).
\end{equation}  

{\bf (ii)}\,lower bound for $R(t)$.
\noindent
The lower bound for $R(t)$ we do not need to separate $R(t)$ into $R_{j}(t)$ with $j = l,h$, and prove at a stroke. The following property is essential.
\begin{equation}\label{412}
\log(1+r^{2}) \leq r^{2}\quad (r \in {\bf R}).
\end{equation}

Once we notice \eqref{412}, the derivation of the lower bound of infinite time blowup rate is similar to \cite{IO-17}. Indeed, by using \eqref{412} and a change of variable one can estimate $R(t)$  as follows.
\begin{align*}
R(t)& = \int_{0}^{\infty}e^{-t\log(1+r^{2})}r^{-1}\sin^{2}(tr)dr\\
& \geq \int_{0}^{\infty}e^{-tr^{2}}r^{-1}\sin^{2}(tr)dr
= \int_{0}^{\infty}e^{-\sigma^{2}}\sigma^{-1}\sin^{2}(\sqrt{t}\sigma)d\sigma.
\end{align*}
By setting 
\[
\nu_{j} := (\frac{1}{4} + j)\frac{\pi}{\sqrt{t}},\qquad \nu_{j}' := (\frac{3}{4} + j)\frac{\pi}{\sqrt{t}}\quad (j = 1,2,3,\cdots)
\]
and using  integration by parts one has
\begin{align}\label{413-2}
R(t) & \geq \frac{1}{2}\sum_{j = 1}^{\infty}\int_{\nu_{j}}^{\nu_{j}'}e^{-\sigma^{2}}\sigma^{-1}d\sigma \geq \frac{1}{4}\int_{\frac{5\pi}{4\sqrt{t}}}^{\infty}e^{-\sigma^{2}}\sigma^{-1}d\sigma \nonumber\\
&=-\frac{1}{4}\log\left(\frac{5\pi}{4\sqrt{t}}\right)e^{-\frac{25\pi^{2}}{16 t}} + \frac{1}{2}\int_{\frac{5\pi}{4\sqrt{t}}}^{\infty}(\sigma\log\sigma)e^{-\sigma^{2}}d\sigma\\
&\geq \frac{1}{8}e^{-\frac{25\pi^{2}}{16 t}}\log t - \frac{1}{4}e^{-\frac{25\pi^{2}}{16 t}}\log \frac{5\pi}{4} - \frac{1}{2}\int_{0}^{\infty}\sigma\vert\log\sigma\vert e^{-\sigma^{2}}d\sigma. \nonumber
\end{align}
\noindent
Since
\[\lim_{t \to \infty}e^{-\frac{25\pi^{2}}{16 t}} = 1,\]
and
\[\int_{0}^{\infty}\sigma\vert\log\sigma\vert e^{-\sigma^{2}}d\sigma < +\infty,\]
\eqref{413-2} implies the desired estimate

\begin{equation}\label{final}
R(t) \geq C\log t\quad (t \gg 1).
\end{equation}

We may note that  \eqref{411} and (4.17) imply the desired statement for $R(t)$.
\hfill
$\Box$

\vspace{0.2cm}
Finally, let us now prove Theorem \ref{main-theo2} at a stroke. \\

{\it Proof of Theorem \ref{main-theo2} completed.}\,It follows from the Plancherel theorem and triangle inequality, with some constant $C_{n} > 0$ one can get
\[C_{n}\Vert u(t,\cdot)\Vert \geq \vert P_{1}\vert\Vert (1+\vert\xi\vert^{2})^{-\frac{t}{2}}\frac{\sin(t\vert\xi\vert)}{\vert\xi\vert}\Vert - \Vert\hat{u}(t,\cdot)-P_{1}(1+\vert\xi\vert^{2})^{-\frac{t}{2}}\frac{\sin(t\vert\xi\vert)}{\vert\xi\vert}\Vert.\]
and
\[C_{n}\Vert u(t,\cdot)\Vert \leq \vert P_{1}\vert\Vert (1+\vert\xi\vert^{2})^{-\frac{t}{2}}\frac{\sin(t\vert\xi\vert)}{\vert\xi\vert}\Vert + \Vert\hat{u}(t,\cdot)-P_{1}(1+\vert\xi\vert^{2})^{-\frac{t}{2}}\frac{\sin(t\vert\xi\vert)}{\vert\xi\vert}\Vert.\]
These inequalities together with Theorem \ref{main-theo}, Lemma \ref{2.8} and Propositions \ref{proposition4.1} and \ref{proposition4.2} imply the desired estimates. This part is, nowadays, well-known (see \cite{I-14, IO-17}). 
\hfill
$\Box$

\par
\vspace{0.5cm}
\noindent{\em Acknowledgement.}
\smallskip
The work of the first author (R. C. CHAR\~AO) was partially supported by PRINT/CAPES - Process 88881.310536/2018-00 and the work of the second author (R. IKEHATA) was supported in part by Grant-in-Aid for Scientific Research (C)15K04958  of JSPS. 



\section{Appendix}

{\small In this appendix, let us describe the outline of proof of the unique existence of a mild solution to problem (1.1)-(1.2) more in detail by applying the Lumer-Phillips Theorem (cf. Pazy \cite[Theorem 4.3]{P}). 

Concerning a relation between two nonnegative self-adjoint operators $A$ and $L$, it holds that $D(A) \subset D(A^{1/2}) \subset D(L) \subset H := L^{2}({\bf R}^{n})$ with  $D(A^{1/2}) = H^{1}({\bf R}^{n})$. We first prepare the Kato-Rellich Theorem.
\begin{theo}{\rm (Kato-Rellich)}\,Let $X$ be a Hilbert space with its norm $\Vert\cdot\Vert$, and let $T: D(T) \subset X \to X$ be a self-adjoint operator in $X$. Furthermore, let $V: D(V) \subset X \to X$ be a symmetric operator in $X$. Assume that\\
{\rm (1)}\,$D(T) \subset D(V)$, \\
{\rm (2)}\,there exist constants $\delta \in [0,1)$ and $C > 0$ such that $\Vert Vu\Vert \leq \delta\Vert Tu\Vert + C\Vert u\Vert$ for $u \in D(T)$.\\
Then, the operator $T + V$ is also self-adjoint in $X$ with its domain $D(T+V) = D(T)$. 
\end{theo}

Now, let ${\cal H}_{0} := H^{1}({\bf R}^{n})\times L^{2}({\bf R}^{n})$ be the Hilbert space with its inner product defined by
\[{\bf <}{u\brack v}, {w\brack z}{\bf >} := (u,w) + (A^{1/2}u,A^{1/2}w) + (v,z),\]
where $(\cdot,\cdot)$ implies the usual inner product in $L^{2}({\bf R}^{n})$. Furthermore, let us define a operator
\[{\cal A}: {\cal H}_{0} \to {\cal H}_{0}\]
by $D({\cal A}) := H^{2}({\bf R})\times H^{1}({\bf R})$, and for $U := \displaystyle{{u\brack v}} \in D({\cal A})$;
\[{\cal A}U := {v\brack -Au-Lv}.\]
Note that $v \in H^{1}({\bf R}^{n})$ implies  $v \in D(L)$. Under these preparations we first show that\\

{\rm (i)}\,The operator ${\cal A}-\frac{1}{2}{\cal I}$ is dissipative in ${\cal H}_{0}$.\\

Indeed, let $U := \displaystyle{{u\brack v}} \in D({\cal A})$. Then  
\[{\bf <}{\cal A}U,U{\bf >} = {\bf <} {v\brack -Au-Lv}, {u\brack v}{\bf >} = (v,u) - (Lv,v) \leq (u,v) \leq \frac{1}{2}(\Vert u\Vert^{2} + \Vert v\Vert^{2}) \leq \frac{1}{2}{\bf <}U,U{\bf >},\]
which implies the desired estimate. Here, one has just used the non-negativity of the self-adjoint operator $L$ in $L^{2}({\bf R}^{n})$ such that 
\[(Lv,v) = c_{n}(\log(1+\vert\cdot\vert^{2})\hat{v},\hat{v}) \geq 0\]
with some constant $c_{n} > 0$.\\

{\rm (ii)}\,For ${\cal B} := {\cal A}-\frac{1}{2}{\cal I}$, we have to check ${\cal R}(\frac{1}{2}{\cal I}-{\cal B}) = {\cal H}_{0}$.\\

Once {\rm (i)} and {\rm (ii)} can be proved, it follows from the Lumer-Phillips Theorem that the operator ${\cal B}$ generates a $C_{0}$ semigroup $e^{t{\cal B}}$ of contractions on ${\cal H}_{0}$, and so $e^{t{\cal A}} = e^{\frac{t}{2}}e^{t{\cal B}}$ can be a generated $C_{0}$ semigroup on ${\cal H}_{0}$ (cf. \cite[Proposition 2.1]{ITY}). 

It suffices to check that ${\cal R}(\frac{1}{2}{\cal I} - {\cal B}) = {\cal R}({\cal I} - {\cal A}) = {\cal H}_{0}$, that is, we have to solve the problem that for each $\displaystyle{{f\brack g}} \in {\cal H}_{0}$, there exists a solution $\displaystyle{{u\brack v}} \in D({\cal A})$ such that
\begin{equation}\label{51}
u-v = f \in H^{1}({\bf R}^{n}),
\end{equation}
\begin{equation}\label{52}
v + Au + Lv = g \in L^{2}({\bf R}^{n}).
\end{equation}
We can find a pair of solution $[u,v]$ to problems \eqref{51} and \eqref{52} by
\begin{equation}\label{53}
u := (I + A + L)^{-1}(f + Lf + g),
\end{equation}
and
\[v := u - f .\]
We easily see that $u \in H^{2}({\bf R}^{n})$, and $v \in H^{1}({\bf R}^{n})$. In order to check the well-posedness of the solution \eqref{53} it is enough to make sure that
 the operator $A + L: D(A) \to L^{2}({\bf R})$ is self-adjoint in $L^{2}({\bf R}^{n})$. 

Let us apply the Kato-Rellich Theorem to check that $A + L$ is self-adjoint in $H := L^{2}({\bf R}^{n})$ with its domain $D(A+L) = D(A)$.

Set $\phi(x) := \log(1+x)(1+x)^{-1}$. Then, since 
\[\max_{x \geq 0}\phi(x) = \phi(e-1) = e^{-1},\]
it holds that for $v \in D(A) = H^{2}({\bf R}^{n})$,
\[
\Vert Lv\Vert^{2} \leq \frac{1}{e^{2}}\int_{{\bf R}^{n}}\vert\hat{v}(\xi)\vert^{2}(1+\vert\xi\vert^{2})^{2}d\xi \leq \frac{4}{e^{2}}(\Vert v\Vert^{2} + \Vert Av\Vert^{2}),
\]
which implies
\begin{equation}
\Vert Lv\Vert \leq \frac{2}{e}(\Vert v\Vert + \Vert Av\Vert)
\end{equation}
with $2/e \in (0,1)$. Therefore, by the Kato-Rellich theorem the operator $A + L$ becomes self-adjoint, and non-negative in $H$. 

Finally, 
\[U(t) = [u(t),u'(t)] := e^{t{\cal A}}[u_{0},u_{1}]\]
becomes a unique mild solution to problem
\[\frac{dU}{dt} = {\cal A}U(t),\quad U(0) = {u_{0}\brack u_{1}}.\]
This implies that the problem (1.1)-(1.2) has a desired unique weak solution 
\[u \in C([0,\infty);H^{1}({\bf R}^{n})) \cap C^{1}([0,\infty);L^{2}({\bf R}^{n})).\]
Finally, by density argument and the multiplier method one can get the energy inequality
\[E_{u}(t) \leq E_{u}(s),\quad (0 \leq s \leq t).\]
}

\end{document}